# Slant Ruled Surfaces

**Mehmet Önder**
*"Independent Researcher" Delibekirli Village, Tepe Street, No. 63, 31440 Kırıkhan, Hatay, Turkey.*
E-mail: mehmetonder197999@gmail.com

**Abstract**
In this study, we define some new types of ruled surfaces called slant ruled surfaces. We give some characterizations for a regular ruled surface to be a slant ruled surface in Euclidean 3-space $E^3$. We show that if the slant ruled surface is developable then the striction curve is a general helix or a slant helix according to the kind of surface. Moreover, we give the relationships between slant ruled surfaces and some offset surfaces such as Bertrand offsets and Mannheim offsets.

**MSC:** 53A25, 14J26.

**Key words:** Ruled surface; Frenet frame; slant ruled surface.

**1. Introduction**

In the local differential geometry, the curves for which the curvatures satisfy some special conditions have an important role. The well-known of such curves are helices. A general helix in Euclidean 3-space $E^3$ is defined by the property that along the curve, the tangent of curve makes a constant angle with a fixed straight line called the axis of the general helix [4]. Therefore, a general helix can be equivalently defined as one whose tangent indicatrix is a planar curve. Certainly, the helices in $E^n$ correspond with those whose unit tangent indicatrices are contained in hyperplanes. A classical result for the helices stated by M.A. Lancret in 1802 and first proved by B. de Saint Venant in 1845 (see [18] for details) is: *A necessary and sufficient condition that a curve to be a general helix is that the ratio of the first curvature to the second curvature be constant i.e., $\kappa/\tau$ is constant along the curve, where $\kappa$ and $\tau$ denote the first and second curvatures of the curve, respectively.* Helices have been studied not only in Euclidean spaces but also in Lorentzian spaces by some mathematicians and different characterizations of these curves have been obtained according to the properties of the spaces [4-7,10,16].

Recently, Izumiya and Takeuchi have introduced the concept of slant helix by saying that the normal lines of the curve make a constant angle with a fixed direction and they have given a characterization of slant helix in Euclidean 3-space $E^3$ [8]. Later, Kula and Yaylı have investigated spherical images, the tangent indicatrix and the binormal indicatrix of a slant helix and they have obtained that the spherical images of a slant helix are spherical helices [11]. Moreover, Kula and et al. have also studied slant helices in Euclidean 3-space and given some other characterizations [12]. Monterde has shown that for a curve with constant curvature and non-constant torsion the principal normal vector makes a constant angle with a fixed constant direction, i.e, the curve is a slant helix [13]. Ali has studied the position vectors of slant helices in Euclidean 3-space [3]. Then Ali and Turgut have given the corresponding characterizations for the position vector of a timelike slant helix in Minkowski 3-space $E_1^3$ [2]. Furthermore, Ali and Lopez have given some new characterizations of slant helices in Minkowski 3-space $E_1^3$ [1].

Moreover, it is well know that a ruled surface has an orthonormal base along its striction line. This frame is called Frenet frame of the ruled surface. In this paper, by considering the Frenet vectors of a ruled surface we give the definitions of some special ruled surfaces for which the Frenet vectors make a constant angle with some fixed directions in the space and we called these surfaces as slant ruled surfaces.

**2. Ruled Surfaces in Euclidean 3-space $E^3$**

In this section, we give a brief summary of ruled surfaces in $E^3$. For more details we refer the ref. [9].

Let $I$ be an open interval in the real line $\mathbb{R}$, $\vec{f} = \vec{f}(u)$ be a regular curve in $E^3$ defined on $I$ and $\vec{q} = \vec{q}(u)$ be a unit direction vector of an oriented line in $E^3$. Then the parametric representation of a ruled surface $N$ is given as follows

$$\vec{r}(u,v) = \vec{f}(u) + v\vec{q}(u). \tag{1}$$

The curve $\vec{f} = \vec{f}(u)$ is called base curve or generating curve of the surface and various positions of the generating lines $\vec{q} = \vec{q}(u)$ are called rulings. In particular, if the direction of $\vec{q}$ is constant, then the ruled surface is said to be cylindrical, and non-cylindrical otherwise.

The distribution parameter of $N$ is defined by

$$d = \frac{\left|\dot{\vec{f}}, \vec{q}, \dot{\vec{q}}\right|}{\left\langle \dot{\vec{q}}, \dot{\vec{q}} \right\rangle}, \tag{2}$$

where $\dot{\vec{f}} = \frac{d\vec{f}}{du}$, $\dot{\vec{q}} = \frac{d\vec{q}}{du}$. If $\left|\dot{\vec{f}}, \vec{q}, \dot{\vec{q}}\right| = 0$, then the normal vectors are collinear at all points of same ruling and at nonsingular points of the surface $N$, the tangent planes are identical, i.e, tangent plane contacts the surface along a ruling. Such a ruling is called a torsal ruling. If $\left|\dot{\vec{f}}, \vec{q}, \dot{\vec{q}}\right| \neq 0$, then the tangent planes of the surface $N$ are distinct at all points of same ruling which is called nontorsal [9].

**Definition 2.1. ([9])** A ruled surface whose all rulings are torsal is called a developable ruled surface. The remaining ruled surfaces are called skew ruled surfaces. From (2) it is clear that a ruled surface is developable if and only if at all its points the distribution parameter is zero.

For the unit normal vector $\vec{m}$ of a ruled surface $N$ we have

$$\vec{m} = \frac{\vec{r}_u \times \vec{r}_v}{\|\vec{r}_u \times \vec{r}_v\|} = \frac{(\dot{\vec{f}} + v\dot{\vec{q}}) \times \vec{q}}{\sqrt{\left\langle \dot{\vec{f}} + v\dot{\vec{q}}, \dot{\vec{f}} + v\dot{\vec{q}} \right\rangle - \left\langle \dot{\vec{f}}, \vec{q} \right\rangle^2}}. \tag{3}$$

Along a ruling $u = u_1$, the unit normal of the surface approaches a limiting direction as $v$ infinitely decreases. This direction is called the asymptotic normal (central tangent) direction and from (3) defined by

$$\vec{a} = \lim_{v \to \pm\infty} \vec{m}(u_1, v) = \frac{\vec{q} \times \dot{\vec{q}}}{\|\dot{\vec{q}}\|}.$$

The point at which $\vec{m}$ is perpendicular to $\vec{a}$ is called the striction point (or central point) and denoted by $C$. The set of striction points of all rulings is called striction curve of the surface. The parametrization of the striction curve $\vec{c} = \vec{c}(u)$ on a ruled surface is given by

$$\vec{c}(u) = \vec{f}(u) + v_0 \vec{q}(u) = \vec{f} - \frac{\left\langle \dot{\vec{q}}, \dot{\vec{f}} \right\rangle}{\left\langle \dot{\vec{q}}, \dot{\vec{q}} \right\rangle} \vec{q}, \tag{4}$$

where $v_0 = -\dfrac{\langle \dot{\vec{q}}, \dot{\vec{f}} \rangle}{\langle \dot{\vec{q}}, \dot{\vec{q}} \rangle}$ is called strictional distance.

The vector $\vec{h}$ defined by $\vec{h} = \vec{a} \times \vec{q}$ is called central normal vector which is the surface normal along the striction curve. Then the orthonormal system $\{C; \vec{q}, \vec{h}, \vec{a}\}$ is called Frenet frame of the ruled surface $N$ where $C$ is the central point and $\vec{q}$, $\vec{h}$, $\vec{a}$ are unit vectors of ruling, central normal and central tangent, respectively.

For the Frenet formulae of the ruled surface $N$ with respect to the arc length $s$ of striction curve we have

$$\begin{bmatrix} d\vec{q}/ds \\ d\vec{h}/ds \\ d\vec{a}/ds \end{bmatrix} = \begin{bmatrix} 0 & k_1 & 0 \\ -k_1 & 0 & k_2 \\ 0 & -k_2 & 0 \end{bmatrix} \begin{bmatrix} \vec{q} \\ \vec{h} \\ \vec{a} \end{bmatrix}, \quad (5)$$

where $k_1 = \dfrac{ds_1}{ds}$, $k_2 = \dfrac{ds_3}{ds}$ and $s_1$, $s_3$ are the arc lengths of the spherical curves circumscribed by the bound vectors $\vec{q}$ and $\vec{a}$, respectively and ruled surfaces satisfying $k_1 \neq 0$, $k_2 = 0$ are called conoids (For details [9]).

***Theorem 2.1.*** ([15]). *Let the striction curve $\vec{c} = \vec{c}(s)$ of the ruled surface $N$ be unit speed i.e., $s$ is arc length parameter of $\vec{c}(s)$ and let $\vec{c}(s)$ be the base curve of the surface. Then $N$ is developable if and only if the unit tangent of the striction curve is the same with the ruling along the curve.*

## 3. $q$-Slant Ruled Surfaces in $E^3$

In this section, we introduce the definition and characterizations of $q$-slant ruled surfaces in $E^3$. First, we give the following definition.

**Definition 3.1.** Let $N$ be a regular ruled surface in $E^3$ given by the parametrization
$$\vec{r}(s, v) = \vec{c}(s) + v\vec{q}(s), \quad \|\vec{q}(s)\| = 1, \quad (6)$$
where $\vec{c}(s)$ is striction curve of $N$ and $s$ is arc length parameter of $\vec{c}(s)$. Let the Frenet frame and non-zero invariants of $N$ be $\{\vec{q}, \vec{h}, \vec{a}\}$ and $k_1, k_2$, respectively. Then, $N$ is called a $q$-slant ruled surface if the ruling makes a constant angle $\theta$ with a fixed non-zero direction $\vec{u}$ in the space, i.e.,
$$\langle \vec{q}, \vec{u} \rangle = \cos\theta = constant; \quad \theta \neq \dfrac{\pi}{2}. \quad (7)$$

Then we give the following characterizations for $q$-slant ruled surfaces. Whenever we talk about $N$, we will mean that the surface has the parametrization and Frenet elements as assumed in Definition 3.1.

***Theorem 3.1.*** *$N$ is a $q$-slant ruled surface if and only if the function*
$$\dfrac{k_1}{k_2} = \tan\theta, \quad (8)$$

*is constant, where $\theta$ is the angle between the ruling $\vec{q}$ and a fixed direction.*

**Proof:** Let $N$ be a $q$-slant ruled surface in $E^3$, $\vec{u}$ be a unit vector of a fixed direction and $\theta$ be constant angle between $\vec{q}$ and $\vec{u}$. Then, $N$ satisfies

$$\langle \vec{q}, \vec{u} \rangle = \cos\theta = constant. \qquad (9)$$

Differentiating (9) with respect to $s$ gives $\langle \vec{h}, \vec{u} \rangle = 0$. Therefore, $\vec{u}$ lies on the plane spanned by the vectors $\vec{q}$ and $\vec{a}$, i.e.,

$$\vec{u} = (\cos\theta)\vec{q} + (\sin\theta)\vec{a}. \qquad (10)$$

By differentiating (10) with respect to $s$ it follows

$$0 = (\cos\theta k_1 - \sin\theta k_2)\vec{h} \qquad (11)$$

and then we have $k_1/k_2 = \tan\theta$ is constant.

Conversely, for a given a regular ruled surface $N$, let the equation (8) is satisfied. We define

$$\vec{u} = (\cos\theta)\vec{q} + (\sin\theta)\vec{a}. \qquad (12)$$

Differentiating (12) and using (8) it follows $\vec{u}' = 0$, i.e., $\vec{u}$ is a constant vector. On the other hand $\langle \vec{q}, \vec{u} \rangle = \cos\theta = constant$. Then $N$ is a $q$-slant ruled surface.

**Theorem 3.2.** $N$ is a $q$-slant ruled surface if and only if $\det(\vec{q}', \vec{q}'', \vec{q}''') = 0$.

**Proof:** Let $N$ be a regular ruled surface in $E^3$. From the Frenet formulae in (5) we have

$$\vec{q}' = k_1 \vec{h},$$
$$\vec{q}'' = -k_1^2 \vec{q} + k_1' \vec{h} + k_1 k_2 \vec{a},$$
$$\vec{q}''' = (-3k_1 k_1')\vec{q} + (k_1'' - k_1^3 - k_1 k_2^2)\vec{h} + (2k_1' k_2 + k_2' k_1)\vec{a},$$

and then

$$\det(\vec{q}', \vec{q}'', \vec{q}''') = k_1^3 k_2^2 \left(\frac{k_1}{k_2}\right)'. \qquad (13)$$

Let now $N$ be a $q$-slant ruled surface in $E^3$. By Theorem 3.1 we have $\dfrac{k_1}{k_2}$ is constant. Then from (13) it follows that $\det(\vec{q}', \vec{q}'', \vec{q}''') = 0$.

Conversely, if $\det(\vec{q}', \vec{q}'', \vec{q}''') = 0$, since the curvatures are non-zero from (13) it is obtained that $\dfrac{k_1}{k_2}$ is constant and Theorem 3.1 gives that $N$ is a $q$-slant ruled surface in $E^3$.

**Theorem 3.3.** $N$ is a $q$-slant ruled surface if and only if $\det(\vec{a}', \vec{a}'', \vec{a}''') = 0$.

**Proof:** From the Frenet formulae we have

$$\vec{a}' = -k_2 \vec{h},$$
$$\vec{a}'' = k_1 k_2 \vec{q} - k_2' \vec{h} - k_2^2 \vec{a},$$
$$\vec{a}''' = (k_1' k_2 + 2k_1 k_2')\vec{q} + (-k_2'' + k_2^3 + k_1^2 k_2)\vec{h} - 3k_2 k_2' \vec{a},$$

and so,

$$\det(\vec{a}', \vec{a}'', \vec{a}''') = k_2^5 \left(\frac{k_1}{k_2}\right)'. \qquad (14)$$

Let now $N$ be a $q$-slant ruled surface in $E^3$. By Theorem 3.1 we have $\dfrac{k_1}{k_2}$ is constant. Then from (14) it follows that $\det(\vec{a}',\vec{a}'',\vec{a}''')=0$.

Conversely, if $\det(\vec{a}',\vec{a}'',\vec{a}''')=0$, since the curvature $k_2$ is non-zero from (14) it is obtained that $\dfrac{k_1}{k_2}$ is constant and Theorem 3.1 gives that $N$ is a $q$-slant ruled surface in $E^3$.

**Theorem 3.4.** $N$ is a $q$-slant ruled surface if and only if
$$\vec{q}''' = m\vec{q}' + 3k_1'\vec{h}', \tag{15}$$
where $m = \dfrac{k_1''}{k_1} - (k_1^2 + k_2^2)$.

**Proof:** Assume that $N$ is a $q$-slant ruled surface. From (5) we get
$$\vec{q}'' = -k_1^2\vec{q} + k_1'\vec{h} + k_1 k_2 \vec{a}, \tag{16}$$
$$\vec{q}''' = (-3k_1 k_1')\vec{q} + (k_1'' - k_1 k_2^2)\vec{h} + (2k_1' k_2 + k_1 k_2')\vec{a} - k_1^2\vec{q}'. \tag{17}$$

Since $\dfrac{k_1}{k_2}$ is constant, by differentiation we have
$$k_1 k_2' = k_2 k_1', \tag{18}$$
and from (5)
$$\vec{h} = \dfrac{1}{k_1}\vec{q}'. \tag{19}$$

Substituting (18) and (19) in (17) gives
$$\vec{q}''' = \left(\dfrac{k_1''}{k_1} - k_1^2 - k_2^2\right)\vec{q}' - (3k_1 k_1')\vec{q} + (3k_2 k_1')\vec{a}. \tag{20}$$

Using the second equation of (5), (15) is obtained from (20).

Conversely, let us assume that (15) holds. Differentiating (19) we obtain
$$\vec{h}' = -\left(\dfrac{k_1'}{k_1^2}\right)\vec{q}' + \left(\dfrac{1}{k_1}\right)\vec{q}'', \tag{21}$$
and so,
$$\vec{h}'' = -\left(\dfrac{k_1'}{k_1^2}\right)'\vec{q}' - 2\left(\dfrac{k_1'}{k_1^2}\right)\vec{q}'' + \left(\dfrac{1}{k_1}\right)\vec{q}'''. \tag{22}$$

Substituting (15) in (22) it follows
$$\vec{h}'' = -2\left(\dfrac{k_1'}{k_1^2}\right)\vec{q}'' - \left[\left(\dfrac{k_1'}{k_1^2}\right)' + \dfrac{m}{k_1}\right]\vec{q}' + 3\left(\dfrac{k_1'}{k_1}\right)\vec{h}'. \tag{23}$$

Now, writing (16) in (23) and using (5) we have
$$\vec{h}'' = -\left[\left(\dfrac{k_1'}{k_1^2}\right)' + \dfrac{m}{k_1}\right]\vec{q}' - k_1'\vec{q} - 2\left(\dfrac{k_1'}{k_1}\right)^2\vec{h} + \left(\dfrac{k_2 k_1'}{k_1}\right)\vec{a}. \tag{24}$$

On the other hand, from (5) it is obtained
$$\vec{h}'' = -k_1\vec{q}' - k_1'\vec{q} - k_2^2\vec{h} + k_2'\vec{a}. \tag{25}$$

Substituting (25) in (24) we have

$$\frac{k_2'}{k_2} = \frac{k_1'}{k_1}. \tag{26}$$

Integrating (26) we get that $\frac{k_1}{k_2}$ is constant and by Theorem 3.1, $N$ is a $q$-slant ruled surface.

***Theorem 3.5.*** *Let $N$ be a developable ruled surface in $E^3$. Then $N$ is a $q$-slant ruled surface if and only if the striction line $\vec{c}(s)$ is a general helix in $E^3$.*

**Proof:** Since $N$ is a developable ruled surface in $E^3$, from Theorem 2.1 we have $\vec{c}'(s) = \vec{t}(s) = \vec{q}(s)$. Then from Definition 3.1, it is clear that $N$ is a $q$-slant ruled surface if and only if the striction line $\vec{c}(s)$ is a general helix in $E^3$.

## 4. $h$-Slant Ruled Surfaces in $E^3$

In this section, we introduce the definition and characterizations of $h$-slant ruled surfaces in $E^3$. First, we give the following definition.

**Definition 4.1.** Let $N$ be a regular ruled surface in $E^3$ given by the parametrization
$$\vec{r}(s,v) = \vec{c}(s) + v\vec{q}(s), \quad \|\vec{q}(s)\| = 1, \tag{27}$$
where $\vec{c}(s)$ is striction curve of $N$ and $s$ is arc length parameter of $\vec{c}(s)$. Let the Frenet frame and non-zero invariants of $N$ be $\{\vec{q}, \vec{h}, \vec{a}\}$ and $k_1, k_2$, respectively. Then $N$ is called a $h$-slant ruled surface if the central normal vector $\vec{h}$ makes a constant angle $\varphi$ with a fixed non-zero direction $\vec{u}$ in the space, i.e.,
$$\langle \vec{h}, \vec{u} \rangle = \cos\varphi = constant; \quad \varphi \neq \frac{\pi}{2}. \tag{28}$$

Then, under the assumptions given in Definition 4.1, we can give the following theorems characterizing $h$-slant ruled surfaces.

***Theorem 4.1.*** *$N$ is a $h$-slant ruled surface if and only if the function*
$$\frac{k_1^2}{\left(k_1^2 + k_2^2\right)^{\frac{3}{2}}} \left(\frac{k_2}{k_1}\right)' \tag{29}$$
*is constant.*

**Proof:** Assume that $N$ is a $h$-slant ruled surface in $E^3$. Let $\vec{u}$ be a fixed constant vector such that $\langle \vec{h}, \vec{u} \rangle = \cos\varphi = c = constant$ where $\varphi$ is the constant angle between $\vec{h}$ and $\vec{u}$. Then for the vector $\vec{u}$ we have
$$\vec{u} = b_1(s)\vec{q}(s) + c\vec{h}(s) + b_2(s)\vec{a}(s), \tag{30}$$
where $b_1 = b_1(s)$ and $b_2 = b_2(s)$ are smooth functions of arc length parameter $s$. Since $\vec{u}$ is constant, differentiation of (30) gives
$$\begin{cases} b_1' - ck_1 = 0, \\ b_1 k_1 - b_2 k_2 = 0, \\ b_2' + ck_2 = 0. \end{cases} \tag{31}$$
From the second equation of system (31) we have

$$b_1 = b_2 \frac{k_2}{k_1}. \tag{32}$$

Moreover,
$$\langle \vec{u}, \vec{u} \rangle = b_1^2 + c^2 + b_2^2 = constant. \tag{33}$$

Substituting (32) in (33) gives
$$b_2^2 \left(1 + \left(\frac{k_2}{k_1}\right)^2\right) = n^2 = constant. \tag{34}$$

If $n = 0$, then $b_2 = 0$ and from (31) we have $b_1 = 0$, $c = 0$. This means that $\vec{u} = \vec{0}$ which is a contradiction. Thus, $n \neq 0$. Then from (34) it is obtained that

$$b_2 = \pm \frac{n}{\sqrt{1 + \left(\frac{k_2}{k_1}\right)^2}}. \tag{35}$$

Considering the third equation of system (31), from (35) we have

$$\frac{d}{ds}\left[\pm \frac{n}{\sqrt{1 + \left(\frac{k_2}{k_1}\right)^2}}\right] = -ck_2. \tag{36}$$

This can be written as

$$\frac{k_1^2}{\left(k_1^2 + k_2^2\right)^{\frac{3}{2}}} \left(\frac{k_2}{k_1}\right)' = \frac{c}{n} = \ell = constant,$$

which is desired.

Conversely, assume that the function in (29) is constant, i.e.,

$$\frac{k_1^2}{\left(k_1^2 + k_2^2\right)^{\frac{3}{2}}} \left(\frac{k_2}{k_1}\right)' = \ell = constant.$$

We define
$$\vec{u} = \frac{k_2}{\sqrt{k_1^2 + k_2^2}} \vec{q} + \ell \vec{h} + \frac{k_1}{\sqrt{k_1^2 + k_2^2}} \vec{a}. \tag{37}$$

Differentiating (37) with respect to $s$ and using (29) we have $\vec{u}' = 0$, i.e., $\vec{u}$ is a constant vector. On the other hand $\langle \vec{h}, \vec{u} \rangle = constant$. Thus, $N$ is a $h$-slant ruled surface in $E^3$.

**Theorem 4.2.** *Let $N$ be a regular ruled surface in $E^3$ with first curvature $k_1 \equiv 1$. Then $N$ is a $h$-slant ruled surface if and only if the second curvature is given by*

$$k_2(s) = \pm \frac{s}{\sqrt{\tan^2 \varphi - s^2}}. \tag{38}$$

**Proof:** Let $N$ be a $h$-slant ruled surface with $k_1 \equiv 1$. Then for a fixed constant unit vector $\vec{u}$ we have
$$\langle \vec{h}, \vec{u} \rangle = \cos \varphi = constant. \tag{39}$$

Differentiating (39) with respect to $s$ gives
$$\langle -\vec{q} + k_2 \vec{a}, \vec{u} \rangle = 0, \qquad (40)$$
and from (40) we have
$$\langle \vec{q}, \vec{u} \rangle = k_2 \langle \vec{a}, \vec{u} \rangle. \qquad (41)$$
If we put $\langle \vec{a}, \vec{u} \rangle = x$, we can write
$$\vec{u} = k_2 x \vec{q} + \cos\varphi \vec{h} + x \vec{a}. \qquad (42)$$
Since $\vec{u}$ is unit, from (42) we have
$$x = \pm \frac{\sin\varphi}{\sqrt{1+k_2^2}}. \qquad (43)$$
Then, the vector $\vec{u}$ is given as follows
$$\vec{u} = \pm \frac{k_2 \sin\varphi}{\sqrt{1+k_2^2}} \vec{q} + \cos\varphi \vec{h} \pm \frac{\sin\varphi}{\sqrt{1+k_2^2}} \vec{a}. \qquad (44)$$
Differentiating (40) with respect to $s$, it follows
$$\langle -(1+k_2^2)\vec{h} + k_2' \vec{a}, \vec{u} \rangle = 0. \qquad (45)$$
Writing $x$ and (39) in (45) we have
$$x = \frac{(1+k_2^2)\cos\varphi}{k_2'}. \qquad (46)$$
From (43) and (46) we obtain the following differential equation,
$$\pm \tan\varphi \frac{k_2'}{\left(1+k_2^2\right)^{3/2}} + 1 = 0. \qquad (47)$$
By integration from (47) we get
$$\pm \tan\varphi \frac{k_2}{\sqrt{1+k_2^2}} + s + c = 0, \qquad (48)$$
where $c$ is integration constant. The integration constant can be subsumed thanks to a parameter change $s \to s - c$. Then (48) can be written as
$$\pm \tan\varphi \frac{k_2}{\sqrt{1+k_2^2}} = -s \qquad (49)$$
which gives us $k_2(s) = \pm \frac{s}{\sqrt{\tan^2\varphi - s^2}}$.

Conversely, assume that $k_2(s) = \pm \frac{s}{\sqrt{\tan^2\varphi - s^2}}$ holds and let us put
$$x = \mp \frac{\sin\varphi}{\sqrt{1+k_2^2}} = \mp \frac{\sin\varphi}{\sqrt{1+\frac{s^2}{\tan^2\varphi - s^2}}} = \mp \cos\varphi \sqrt{\tan^2\varphi - s^2}, \qquad (50)$$
where we are assuming that when $k_2$ has the positive (negative) sign, then $x$ gets the negative (positive) sign and $\theta$ is constant. Thus, $k_2 x = -s\cos\varphi$. Let now consider the vector $\vec{u}$ defined by
$$\vec{u} = \cos\varphi \left( s\vec{q} + \vec{h} \pm \left( \sqrt{\tan^2\varphi - s^2} \right) \vec{a} \right) \qquad (51)$$

We will prove that $\vec{u}$ is constant and makes a constant angle $\varphi$ with $\vec{h}$. By differentiating (51) and using Frenet formulae we have $\vec{u}' = 0$, i.e., the direction of $\vec{u}$ is constant and $\langle \vec{h}, \vec{u} \rangle = \cos\varphi = constant$. Then $N$ is a $h$-slant ruled surface.

On the other hand, if the striction line $\vec{c}(s)$ is a geodesic on $N$, then the principal normal vector $\vec{n}$ of $\vec{c}(s)$ and the central normal vector $\vec{h}$ of $N$ coincide. Then, we have the following corollary.

***Corollary 4.1.*** *Let the striction line $\vec{c}(s)$ be a geodesic on $N$. Then $N$ is a $h$-slant ruled surface if and only if the striction line is a slant helix in $E^3$.*

If the ruled surface $N$ is developable, then by Theorem 2.1, the Frenet frame $\{\vec{t}, \vec{n}, \vec{b}\}$ of the striction line $\vec{c}(s)$ coincides with the frame $\{\vec{q}, \vec{h}, \vec{a}\}$ and we can give the following corollary.

***Corollary 4.2.*** *Let $N$ be a developable surface. Then $N$ is a $h$-slant ruled surface if and only if the striction line is a slant helix in $E^3$.*

## 5. $a$-Slant Ruled Surfaces in $E^3$

In this section we introduce the definition of $a$-slant ruled surfaces in $E^3$.

**Definition 5.1.** Let $N$ be a regular ruled surface in $E^3$ given by the parametrization
$$\vec{r}(s,v) = \vec{c}(s) + v\vec{q}(s), \quad \|\vec{q}(s)\| = 1,$$
where $\vec{c}(s)$ is striction curve of $N$ and $s$ is arc length parameter of $\vec{c}(s)$. Let the Frenet frame and non-zero invariants of $N$ be $\{\vec{q}, \vec{h}, \vec{a}\}$ and $k_1, k_2$, respectively. Then $N$ is called a $a$-slant ruled surface if the central tangent vector $\vec{a}$ makes a constant angle $\mu$ with a fixed non-zero direction $\vec{u}$ in the space, i.e.,
$$\langle \vec{a}, \vec{u} \rangle = \cos\mu = constant; \quad \mu \neq \frac{\pi}{2}.$$

From (10) it is celar that a ruled surface $N$ is $a$-slant ruled surface if and only if it is a $q$-slant ruled surface. So, all the theorems given in Section 3 also characterize the $a$-slant ruled surfaces.

After these definitions and characterizations of slant ruled surfaces we have the followings:

Let $N_1$ and $N_2$ be two ruled surfaces in $E^3$ with Frenet frames $\{\vec{q}_1, \vec{h}_1, \vec{a}_1\}$ and $\{\vec{q}_2, \vec{h}_2, \vec{a}_2\}$, respectively. If $N_1$ and $N_2$ have common central normals i.e., $\vec{h}_1 = \vec{h}_2$ at the corresponding points of their striction lines, then $N_1$ and $N_2$ are called Bertrand offsets [17]. Similarly, if $\vec{a}_1 = \vec{h}_2$ at the corresponding points of their striction lines, then the surface $N_2$ is called a Mannheim offset of $N_1$ and the ruled surfaces $N_1$ and $N_2$ are called Mannheim offsets [14]. Considering these definitions we come to the following corollaries:

***Corollary 5.1.*** *Let $N_1$ be a $h$-slant ruled surface. Then the Bertrand offsets of $N_1$ form a family of $h$-slant ruled surfaces.*

***Corollary 5.2.*** *Let $N_1$ and $N_2$ form a Mannheim offset. Then $N_1$ is a $q$-slant (or $a$-slant) ruled surface if and only if $N_2$ is a $h$-slant ruled surface.*

## 6. Examples

**Example 6.1.** Let consider the ruled surface $N$ given by the parametrization

$$r(s,v) = \left(\frac{1}{3}(1+s)^{3/2} + v\frac{1}{2}(1+s)^{1/2}, \frac{1}{3}(1-s)^{3/2} - v\frac{1}{2}(1-s)^{1/2}, \frac{1}{\sqrt{2}}s + v\frac{1}{\sqrt{2}}\right).$$

(Fig. 6.1). It is easily seen that $\det(\vec{q}',\vec{q}'',\vec{q}''') = 0$. Then Theorem 3.2 gives that $N$ is a $q$-slant ruled surface in $E^3$. Since the $q$-slant ruled surfaces are also $a$-slant, $N$ is also a $a$-slant ruled surface.

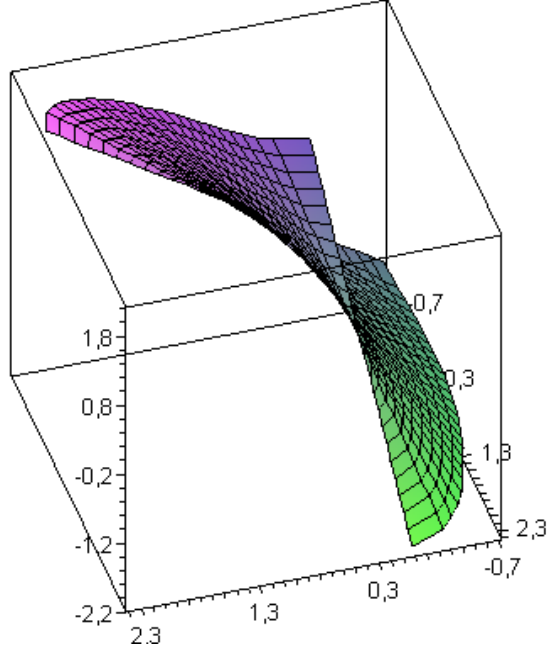

**Fig. 6.1.**

**Example 6.2.** Let consider the ruled surface $N$ given by the parametrization $\vec{r}(s,v) = (r_1, r_2, r_3)$ where

$$r_1(s,v) = \frac{25}{612}\sin(18s) - \frac{9}{1700}\sin(50s) + v\left[\frac{50}{68}\cos(18s) - \frac{18}{68}\cos(50s)\right],$$

$$r_2(s,v) = -\frac{25}{612}\cos(18s) + \frac{9}{1700}\cos(50s) + v\left[\frac{50}{68}\sin(18s) - \frac{18}{68}\sin(50s)\right],$$

$$r_3(s,v) = \frac{15}{272}\sin(16s) + v\frac{15}{17}\cos(16s).$$

(Fig. 6.2). After some computations, the curvatures of the surface are obtained as

$$k_1 = \frac{510}{17}\sin(16s), \quad k_2 = \frac{510}{17}\cos(16s).$$

Then we have

$$\frac{k_1^2}{(k_1^2+k_2^2)^{\frac{3}{2}}}\left(\frac{k_2}{k_1}\right)' = -\frac{136}{255},$$

and Theorem 4.1 gives that $N$ is a $h$-slant ruled surface.

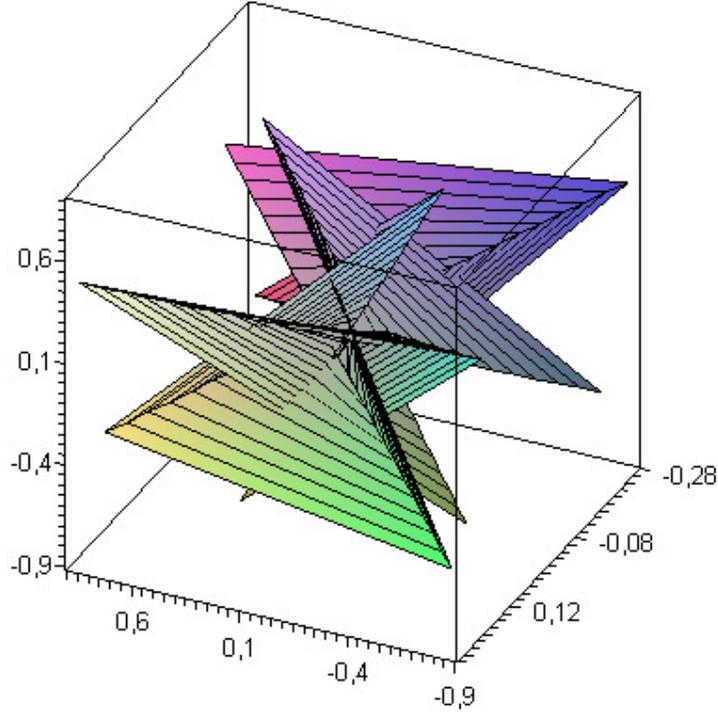

**Fig. 6.2.**

## 7. Conclusions

New types of ruled surfaces in Euclidean 3-space $E^3$ are defined and called slant ruled surfaces. Some properties of these special surfaces are obtained and the relationships between slant ruled surfaces and their striction lines are introduced. Of course, one can consider the definitions and characterizations given in this study for the ruled surfaces of the other spaces and obtain the corresponding characterizations.